\documentclass[12pt]{article}
 \usepackage{amssymb,amsmath,amsfonts}%,numbysec}
 \usepackage{graphics,graphicx,color}
 -\marginparwidth \oddsidemargin -4mm \evensidemargin
 \oddsidemargin%
 \usepackage{amsmath}    % remove % for AMS-LaTeX
 \advance\hoffset by 5mm

 \title{On the asymptotic stability of bound states in 2D cubic Schr\" odinger equation}
 \author{E. Kirr \thanks{Department of Mathematics, University of Illinois at Urbana-Champaign}
 \hspace{0.05in} and A. Zarnescu \thanks{Department of Mathematics, University of Chicago, Chicago, IL}}
 \date{\today}
 \newtheorem{lemma}{Lemma}[section]
 \newtheorem{theorem}{Theorem}[section]
 \newtheorem{proposition}{Proposition}[section]
 \newtheorem{corollary}{Corollary}[section]

\begin{document}

 \maketitle
 \begin{abstract}
 \noindent We consider the cubic nonlinear Schr\"{o}dinger equation in two space
 dimensions with an attractive potential. We study the asymptotic
 stability of the nonlinear bound states, i.e. periodic in time
 localized in space solutions. Our result shows that all solutions
 with small, localized in space initial data, converge to the set of
 bound states. Therefore, the center manifold in this problem is a
 global attractor. The proof hinges on dispersive estimates that we
 obtain for the non-autonomous, non-Hamiltonian, linearized dynamics
 around the bound states.
 \end{abstract}

\section{Introduction}

In this paper we study the long time behavior of solutions of the
cubic nonlinear Schr\" odinger equation (NLS) with potential in two
space dimensions (2-d): {\setlength\arraycolsep{2pt}
\begin{eqnarray}
  i\partial_t u(t,x)&=&[-\Delta_x+V(x)]u+\gamma|u|^2u,\ \ \ \ t>0,\
  x\in\mathbb{R}^2\label{eq:ufull}\\
  u(0,x)&=&u_0(x)\label{eq:ic}
  \end{eqnarray}}
  where $\gamma\in\mathbb{R}-\{0\}.$
The equation has important applications in statistical physics,
optics and water waves. It describes certain limiting behavior of
Bose-Einstein condensates \cite{dgps:bec,lsy:2d} and propagation of
time harmonic waves in wave guides \cite{kn:kp,kn:Marcuse,nm:no}. In
the latter, $t$ plays the role of the coordinate along the axis of
symmetry of the wave guide.

It is well known that this nonlinear equation admits periodic in
time, localized in space solutions (bound states or solitary waves).
They can be obtained via both variational techniques
\cite{bl:i,str:sw,rw:bs} and bifurcation methods \cite{pw:cm,rw:bs},
see also next section. Moreover the set of periodic solutions can be organized as
a manifold (center manifold). Orbital stability of solitary waves, i.e. stability
modulo the group of symmetries $u\mapsto e^{-i\theta}u,$ was first proved in
\cite{rw:bs,mw:ls}, see also
\cite{gss:i,gss:ii,ss:ins}.
%This means that if the
%initial condition is small and close to a periodic solution in the energy space
%$H^1$ then the solution remains close to the set of all periodic solutions.

In this paper we are going to show that the center manifold is in
fact a global attractor for all small, localized in space initial
data. This means that the solution decomposes into a modulation of
periodic solutions (motion on the center manifold) and a part that
decays in time via a dispersion mechanism (radiative part). For a
precise statement of hypotheses and the result see Section
\ref{se:result}.

Asymptotic stability studies of solitary waves were initiated in the
work of A. Soffer and M. I. Weinstein \cite{sw:mc1,sw:mc2}, see also
\cite{bp:asi,bp:asii,bs:as,sc:as,gnt:as}. Center manifold analysis
was introduced in \cite{pw:cm}, see also \cite{kn:Wed}. The
techniques developed in these papers do not apply to our problem.
Indeed the weaker $L^1\rightarrow L^\infty$ dispersion estimates for
Schr\" odinger operators in 2-d, see (\ref{est:Lp}), compared to 3-d
and higher, respectively lack of end point Strichartz estimates in
$d=2,$ prevent the bootstrapping argument in
\cite{sc:as,pw:cm,sw:mc1,sw:mc2}, respectively \cite{gnt:as},  from
closing. The technique of virial theorem, used in
\cite{bp:asi,bp:asii,bs:as} to compensate for the weak dispersion in
1-d, would require at least a quintic nonlinearity in our 2-d case.
Finally, in \cite{kn:Wed}, the nonlinearity is localized in space, a
feature not present in our case, which allows the author to
completely avoid any $L^1\rightarrow L^\infty$ estimates.

To overcome this difficulties we used Strichartz estimates, fixed point and
interpolation techniques to carefully analyze the full, time dependent,
non-Hamiltonian, linearized dynamics around solitary waves. We obtained dispersive
estimates that are similar with the ones for the time independent, Hamiltonian
Schr\" odinger operator, see section \ref{se:lestimates}. Related results have
been proved for the 1-d and 3-d case in \cite{rs:tdde,ks:sm1,ws:sm3} but their
argument does not extend to the 2-d case. We relied on these estimates to
understand the nonlinear dynamics via perturbation techniques. We think that our
estimates are also useful in approaching the dynamics around large 2-d solitary
waves while the techniques that we develop may be used in lowering the power of
nonlinearity needed for the asymptotic stability results in 1 and 3-d mentioned in
the previous paragraph.

Note that, in 3-d, the case of a center manifold formed by two distinct branches
(ground state and excited state) has been analyzed. Under the assumption that the
excited branch is sufficiently far away from the ground state one, in a series of
papers \cite{sw:sgs,ty:ad,ty:rel,ty:sd}, the authors show asymptotic stability of
the ground states with the exception of a finite dimensional manifold where the
solution converges to excited states. We cannot extend such a result to our 2-d
problem as of now. The reason is the slow convergence in time towards the center
manifold, $t^{-1+}$ in 2-d compared to $t^{-3/2}$ in 3-d. This prevents us from
even analyzing the projected dynamics on a single branch center manifold, i.e. the
evolution of one complex parameter describing the projection of the solution on
the center manifold, and obtain, for example, convergence to a periodic orbit as
in
\cite{bs:as, pw:cm, sw:mc2}. However the evolution of this parameter, respectively
two parameters in the presence of the excited branch, is given by an ordinary
differential equation (ODE), respectively a system of two ODE's, and the
contribution of most of the terms can be determined from our estimates, see the
discussion in Section~\ref{se:conclusions}. We think it is only a matter of time
until the remaining ones will be understood.

The paper is organized as follows. In the next section we discuss
previous results regarding the manifold of periodic solutions that
we subsequently need. In section \ref{se:result} we formulate and
prove our main result. As we mentioned before the proof relies on
certain estimates for the linear dynamics which we prove in section
\ref{se:lestimates}. We conclude with possible extensions and
comments in section \ref{se:conclusions}.

\bigskip

\noindent{\bf Notations:} $H=-\Delta+V;$

$L^p=\{f:\mathbb{R}^2\mapsto \mathbb{C}\  |\ f\ {\rm measurable\
and}\ \int_{\mathbb{R}^2}|f(x)|^pdx<\infty\},$
$\|f\|_p=\left(\int_{\mathbb{R}^2}|f(x)|^pdx\right)^{1/p}$ denotes
the standard norm in these spaces;

$<x>=(1+|x|^2)^{1/2},$ and for $\sigma\in\mathbb{R},$ $L^2_\sigma$
denotes the $L^2$ space with weight $<x>^{2\sigma},$ i.e. the space
of functions $f(x)$ such that $<x>^{\sigma}f(x)$ are square
integrable endowed with the norm
$\|f(x)\|_{L^2_\sigma}=\|<x>^{\sigma}f(x)\|_2;$

$\langle f,g\rangle =\int_{\mathbb{R}^2}\overline f(x)g(x)dx$ is the
scalar product in $L^2$ where $\overline f=$ the complex conjugate
of the complex number $f;$

$P_c$ is the projection on the continuous spectrum of $H$ in $L^2;$

$H^n$ denote the Sobolev spaces of measurable functions having all
distributional partial derivatives up to order $n$ in $L^2,
\|\cdot\|_{H^n}$ denotes the standard norm in this spaces.

\section{Preliminaries. The center manifold.}\label{se:prelim}

The center manifold is formed by the collection of periodic solutions for
(\ref{eq:ufull}):
\begin{equation}\label{eq:per}
  u_E(t,x)=e^{-iEt}\psi_E(x)
\end{equation}
where $E\in\mathbb{R}$ and $0\not\equiv\psi_E\in H^2(\mathbb{R}^2)$
satisfy the time independent equation:
\begin{equation}\label{eq:ev}
[-\Delta+V]\psi_E+\gamma|\psi_E|^2\psi_E=E\psi_E
\end{equation}
Clearly the function constantly equal to zero is a solution of (\ref{eq:ev}) but
(iii) in the following hypotheses on the potential $V$ allows for a bifurcation
with a nontrivial, one parameter family of solutions:

\bigskip\noindent{\bf (H1)} Assume that
\begin{itemize}
  \item[(i)] There exists $C>0$ and $\rho >3$ such that:
  $$|V(x)|\le C<x>^{-\rho},\quad {\rm for\ all}\ x\in\mathbb{R}^2;$$
  \item[(ii)] $0$ is a regular point\footnote{see
 \cite[Definition 7]{ws:de2} or $M_\mu=\{0\}$ in relation (3.1) in \cite{mm:ae}} of the
spectrum of
  the linear operator $H=-\Delta+V$ acting on $L^2;$
  \item [(iii)]$H$ acting on $L^2$ has exactly one
  negative eigenvalue $E_0<0$ with corresponding normalized
  eigenvector $\psi_0.$ It is well known that $\psi_0(x)$ can be
  chosen strictly positive and exponentially decaying as
$|x|\rightarrow\infty.$
\end{itemize}

\par\noindent Conditions (i)-(ii) guarantee the applicability of dispersive
estimates of Murata \cite{mm:ae} and Schlag \cite{ws:de2} to the
Schr\" odinger group $e^{-iHt},$ see section \ref{se:lestimates}. In
particular (i) implies the local well posedness in $H^1$ of the
initial value problem (\ref{eq:ufull}-\ref{eq:ic}), see section
\ref{se:result}.

Condition (iii) guarantees bifurcation of nontrivial solutions of \eqref{eq:ev}
from $(E_0,\psi_0).$ In Section \ref{se:conclusions}, we discuss  the possible
effects of relaxing (iii) to allow for finitely many negative eigenvalues. We
construct the center manifold by applying the standard bifurcation argument in
Banach spaces
\cite{ln:fa} for (\ref{eq:ev}) at $E=E_0.$ We follow
\cite{pw:cm} and decompose the solution of (\ref{eq:ev}) in its projection onto
the discrete and continuous part of the spectrum of $H:$
$$\psi_E=a\psi_0+h,\quad a=\langle \psi_0,\psi_E\rangle,\
h=P_c\psi_E.$$ Using the notations {\setlength\arraycolsep{2pt}
\begin{eqnarray}
 f_p(a,h)&\equiv
 &\langle\psi_0,|a\psi_0+h|^2(a\psi_0+h)\rangle, \label{eq:fp}\\
f_c(a,h)&\equiv
 &P_c|a\psi_0+h|^2(a\psi_0+h),\label{eq:fc}
\end{eqnarray}}
and projecting (\ref{eq:ev}) onto $\psi_0$ and its orthogonal
complement $={\rm Range}\ P_c$ we get: {\setlength\arraycolsep{2pt}
\begin{eqnarray}
h&=&-\gamma(H-E)^{-1}f_c(a,h)\label{eq:evc}\\
E_0-E&=&-\gamma a^{-1}f_p(a,h).\label{eq:evp}
\end{eqnarray}}
Although we are using milder hypothesis on $V$ the argument in the Appendix of
\cite{pw:cm} can be easily adapted to show that:
$${\cal F}(E,a,h)=h+\gamma(H-E)^{-1}f_c(a,h)$$
is a $C^1$ function from $(-\infty,0)\times\mathbb{C}\times
L^2_\sigma\cap H^2$ to $L^2_\sigma\cap H^2$ and ${\cal
F}(E_0,0,0)=0,$ $D_h{\cal F}(E_0,0,0)=I.$ Therefore the implicit
function theorem applies to equation (\ref{eq:evc}) and leads to the
existence of $\delta_1>0$ and the $C^1$ function $\tilde h(E,a)$
from $(E_0-\delta_1,E_0+\delta_1)\times \{a\in\mathbb{C}\ :\
|a|<\delta_1\} $ to $L^2_\sigma\cap H^2$ such that (\ref{eq:evc})
has a unique solution $h=\tilde h(E,a)$ for all $E\in
(E_0-\delta_1,E_0+\delta_1)$ and $|a|<\delta_1.$ Note that if
$(a,h)$ solves (\ref{eq:evc}) then $(e^{i\theta}a,e^{i\theta}h),\
\theta\in[0,2\pi )$ is also a solution, hence by uniqueness we have:
\begin{equation}\label{eq:hsym}
\tilde h(E,a)=\frac{a}{|a|}\tilde h(E,|a|).
\end{equation}
Because $\psi_0$ is real valued, we could apply the implicit function theorem to
(\ref{eq:evc}) under the restriction $a\in\mathbb{R}$ and $h$ in the subspace of
real valued functions as it is actually done in \cite{pw:cm}. By uniqueness of the
solution we deduce that $\tilde h(E,|a|)$ is a real valued function.

Replacing now $h=\tilde h(E,a)$ in (\ref{eq:evp}) and using
(\ref{eq:fp}) and (\ref{eq:hsym}) we get the equivalent formulation:
\begin{equation}\label{eq:evp1}
E_0-E=-\gamma |a|^{-1}f_p(|a|,\tilde h(E,|a|)).
\end{equation}
To this we can apply again the implicit function theorem by
observing that $ G(E,a)=E_0-E+\gamma |a|^{-1}f_p(|a|, \tilde
h(E,|a|))$ is a $C^1$ function \cite[Appendix]{pw:cm} from
$(E_0-\delta_1,E_0+\delta_1)\times (-\delta_1,\delta_1)$ to
$\mathbb{R}$ with the properties $G(E_0,0)=0,$ $\partial_E
G(E_0,0)=-1.$ We obtain the existence of $0<\delta\le\delta_1,$
$0<\delta_E\le \delta_1$ and the $C^1$ function $\tilde E:(-\delta ,
\delta)\mapsto (E_0-\delta_E,E_0+\delta_E)$  such that, for
$|E-E_0|<\delta_E,\ |a|<\delta,$ the unique solution of
(\ref{eq:evp}) with $h=\tilde h(E,a),$ is given by $E=\tilde
E(|a|).$ If we now define:
$$h(a)\equiv\frac{a}{|a|}\tilde h(E(|a|,|a|)$$
we have the following center manifold result:

\begin{proposition}\label{pr:cm} There exist $\delta_E,\delta>0$ and
the $C^1$ function
$$h:\{a\in\mathbb{C}\ :\ |a|<\delta\}\mapsto L^2_\sigma\cap H^2
,$$ such that for $|E-E_0|<\delta_E$ the eigenvalue problem
(\ref{eq:ev}) has a unique solution up to multiplication with
$e^{i\theta},\ \theta\in [0,2\pi),$ which can be represented as:
$$\psi_E=a\psi_0+h(a),\qquad \langle\psi_0,h(a)\rangle =0,\
|a|<\delta .$$
\end{proposition}

Since $\psi_0(x)$ is exponentially decaying as $|x|\rightarrow\infty$ the
proposition implies that $\psi_E\in L^2_\sigma .$ A regularity argument, see
\cite{sw:mc1}, gives a stronger result:

\begin{corollary}\label{co:decay} For any $\sigma\in\mathbb{R},$
there exists a finite constant $C_\sigma$ such that:
$$\|<x>^\sigma\psi_E\|_{H^2}\le C_\sigma\|\psi_E\|_{H^2}.$$
\end{corollary}

We are now ready to prove our main result.

\section{Main Result. The collapse on the center
manifold.}\label{se:result}

 \begin{theorem}\label{theorem:main} Assume that hypothesis (H1) is valid and fix $\sigma
>2.$ Then there
 exists an $\varepsilon_0>0$ such that for all initial conditions
 $u_0(x)$ satisfying
 $$\max\{\|u_0\|_{L^2_{\sigma}},\|u_0\|_{H^1}\}\le \varepsilon_0
 $$ the initial value problem
(\ref{eq:ufull})-(\ref{eq:ic}) is
 globally well-posed in $H^1.$

 \par Moreover, for all $t\in\mathbb{R}$ and $p>2,$ we have that:{\setlength\arraycolsep{2pt}
 \begin{eqnarray}
u(t,x)&=&\underbrace{a(t)\psi_0(x)+h(a(t))}_{\psi_{E}(t)}+r(t,x)\label{eq:udec}\\
 \|r(t)\|_{L^2_{-\sigma}}&\le &\frac{\bar
 C_1\varepsilon_0}{(1+|t|)^{1-2/p}}\nonumber\\
 \|r(t)\|_{L^p}&\le &\frac{\bar C_2\varepsilon_0
  \log(2+|t|)}{(1+|t|)^{1-2/p}}\nonumber
  %\\
 %\|r(t)\|_{L^2}&\le &\bar C_3\varepsilon_0\nonumber
  \end{eqnarray}} with the constants
  $\bar C_1,\bar C_2,$%\bar C_3$
  independent of $\varepsilon_0$ and $\bar C_2$ depending on $p>2$.
 \end{theorem}

\smallskip\par Before proving the theorem let us note that
(\ref{eq:udec}) decomposes the evolution of the solution of
(\ref{eq:ufull})-(\ref{eq:ic}) into an evolution on a center manifold $\psi_E(t)$
and the ``distance" from the center manifold $r(t).$ The estimates on the latter
show collapse of solution onto the center manifold. The evolution on the center
manifold is determined by equation (\ref{eq:at}) below. We discuss it in Section
\ref{se:conclusions}.

\smallskip\par{\bf Proof of Theorem \ref{theorem:main}.} It is well
known that under hypothesis (H1)(i) the initial value problem
(\ref{eq:ufull})-(\ref{eq:ic}) is locally well posed in the energy
space $H^1$ and its $L^2$ norm is conserved, see for example \cite[
Corollary 4.3.3. at p. 92]{caz:bk}. Global well posedness follows
via energy estimates from $\|u_0\|_2$ small, see \cite[Remark 6.1.3
at p. 165]{caz:bk}.

In particular we can define
$$a(t)=\langle\psi_0,u(t)\rangle,\quad {\rm for}\ t\in\mathbb{R}.$$
Cauchy-Schwarz inequality implies
$$|a(t)|\le \|u(t)\|_2\|\psi_0\|_2=\|u_0\|_2\le\varepsilon_0,\ {\rm
for}\ t\in\mathbb{R}$$ where we also used conservation of $L^2$ norm of $u.$
Hence, if we choose $\varepsilon_0<\delta$ we can define $h(a(t)),\
t\in\mathbb{R},$ see Proposition \ref{pr:cm}. We then obtain (\ref{eq:udec}) where
$$r(t)=u(t)-a(t)\psi_0-h(a(t)),\quad \langle\psi_0,r(t)\rangle\equiv
0.$$ The solution is now described by the scalar $a(t)\in\mathbb{C}$
and $r(t)\in C(\mathbb{R},H^1)$. Their equations are obtained by
projecting (\ref{eq:ufull}) onto $\psi_0$ and its orthogonal
complement in $L^2:$ {\setlength\arraycolsep{2pt}
 \begin{eqnarray}
 i\frac{da}{dt}&=&E(|a(t)|)a(t)+\gamma\langle\psi_0,2|\psi_E|^2
r+\psi_E^2\bar r+
 2\psi_E|r|^2+\bar\psi_E r^2+|r|^2 r\rangle\label{eq:at}\\
 i\frac{\partial r}{\partial t}&=&Hr+\gamma P_c[2|\psi_E|^2
r+\psi_E^2\bar r
 +2\psi_E|r|^2+\bar\psi_E r^2+|r|^2 r]\nonumber\\
 &&-\gamma Dh|_{a(t)}\langle\psi_0,2|\psi_E|^2 r+\psi_E^2\bar r+
 2\psi_E|r|^2+\bar\psi_E r^2+|r|^2 r\rangle
 \label{eq:r1}
 \end{eqnarray}}
where we used the identities $\psi_E=a\psi_0+h(a)$ and $Dh|_{a}[a]=E(|a|)h(a)$ for
$a\in\mathbb{C},\ |a|<\delta.$

In order to obtain the estimates for $r(t),$ we analyze equation
(\ref{eq:r1}). In the next section we study its linear part:
$$
 \begin{cases}\setlength\arraycolsep{2pt}
 i\frac{\partial z}{\partial t}&=Hz+\gamma P_c[2|\psi_E|^2z+\psi_E^2
\bar z-\gamma Dh|_{a(t)}\langle\psi_0,2|\psi_E|^2 z+\psi_E^2\bar
 z\rangle]\\
 z(s)&=v
 \end{cases}
$$

 \par Let us denote by $\Omega(t,s)v$ the operator which
 associates to the function $v$ the solution of the above equation:
 \begin{equation}
 \Omega(t,s)v\stackrel{\textrm{def}}{=}z
 \label{def:Omega}
 \end{equation}

\noindent The estimates that we need for this linear propagator are
proved in the next section.
\par Now, using Duhamel's principle  (\ref{eq:r1}) becomes
\begin{eqnarray}
r(t)=\Omega(t,0)r(0)+\int_0^t \Omega(t,s)\{\gamma
P_c[2\psi_E|r(s)|^2+\bar\psi_E r^2(s)+|r|^2 r(s)]\nonumber\\
-\gamma Dh|_{a(t)}[\langle\psi_0,2\psi_E|r|^2+\bar\psi_E
r^2+|r|^2r\rangle]\}ds \label{eq:r2}
\end{eqnarray}
It is here where we differ essentially from the approach for the 1-d
case \cite{bp:asi,bp:asii,bs:as} and 3-d case
\cite{sc:as,sw:mc1,sw:mc2,pw:cm}. The right hand side of our
equation contains only nonlinear terms in $r.$ Hence, if we make the
ansatz $r(t)\sim (1+t)^{-3/4}$ then the quadratic and cubic terms in
(\ref{eq:r2}) decay like $(1+s)^{-6/4}$ respectively $
(1+s)^{-9/4}.$ Both are integrable functions in time, hence, via
convolution estimates, the integral term on the right hand side
decays like $\Omega(t,0).$ We have a chance of "closing" the ansatz
provided $\Omega(t,0)\sim (1+t)^{-3/4}.$ Contrast this with the case
in the above cited papers where a linear term in $r$ is present on
the right hand side. Same argument leads to a loss of $1/4$ power
decay in the linear term and requires $\Omega(t,0)\equiv
e^{-iHt}\sim (1+t)^{-1-\delta},\ \delta>0.$ for closing. This turns
out to be impossible in $L^p$ norms in 2-d, see (\ref{est:Lp}),
while the use of weighted $L^2$ norms, see (\ref{Murata}), for
delocalized terms as the cubic term in (\ref{eq:r2}), would require
compensation via virial inequalities, see \cite{bp:asi}, which needs
a much higher power nonlinearity than cubic in the delocalized
terms\footnote{Heuristically we arrived at quintic power
nonlinearity, hence this technique may be applicable to the quintic
Schr\" odinger in 2-D but definitely not to the cubic one.}.

The following Lemma makes the above heuristic argument rigorous:
\begin{lemma}  There exists
$\varepsilon_1>0$ and $\varepsilon_2>0$ such that if
$\|<x>^\sigma\psi_E\|_{H^2}\le \varepsilon_1$ and the initial
condition $r(0)$ satisfies
$$
\max\{\|r(0)\|_{L^2_\sigma},\|r(0)\|_{H^1}\}<\varepsilon_2 $$  the
initial value problem (\ref{eq:r2}) is globally well-posed in
$C(\mathbb{R},L^2_{-\sigma}\cap L^2\cap L^p),\ 6\le p<\infty$ and for all
$t\in\mathbb{R}:$ {\setlength\arraycolsep{2pt}
\begin{eqnarray}
  \|r(t)\|_{L^2_{-\sigma}} &\le& \frac{\bar
C_1\varepsilon_0}{(1+|t|)^{1-2/p}}\nonumber\\
  \|r(t)\|_{L^p} &\le& \frac{\bar C_2\varepsilon_0
\log(2+|t|)}{(1+|t|)^{1-2/p}}\nonumber\\
  \|r(t)\|_{L^2} &\le& \bar C_3\varepsilon_0\nonumber\\
  \nonumber
\end{eqnarray}} with $\bar C_1,\bar C_2,\bar C_3$
  independent of $\varepsilon_0$ and $\bar C_2$ depending on $p\ge 6$.
\label{lemma:r}
\end{lemma}

Note that the Lemma finishes the proof of the theorem. Indeed, we now have two
solutions of (\ref{eq:r1}), one in $C(\mathbb{R},H^1)$ from classical well
posedness theory and one in $C(\mathbb{R},L^2_{-\sigma}\cap L^2\cap L^p),\ p\ge 6$
from the Lemma. Using uniqueness and the continuous embedding of $H^1$ in
$L^2_{-\sigma}\cap L^2\cap L^p,$ we infer that the two solutions must coincide.
Therefore, the time decaying estimates in the Lemma hold also for the $H^1$
solution. The $L^p,\ 2<p<6$ estimates in the theorem follow from interpolation:
$$\|r(t)\|_{L^p}\le\|r(t)\|_{L^2}^{3/p-1/2}\|r(t)\|_{L^6}^{3/2-3/p}\le\frac{\bar C_3\bar C_2\varepsilon_0
\log(2+|t|)}{(1+|t|)^{1-2/p}}$$
It remains to prove the Lemma:

\smallskip\par\noindent {\bf Proof of  Lemma~\ref{lemma:r}.}
\par\noindent Fix $p\ge 6.$ We will show that (\ref{eq:r2}) has a solution by
applying the contraction principle in the functional space
\begin{eqnarray}
 Y=\{f:R\rightarrow L^2_{-\sigma}\cap L^p\cap L^2|\,\,
\sup_{t\ge
0}{(1+|t|)^{1-\frac{2}{p}}}\|f(t)\|_{L^2_{-\sigma}}<\infty
\nonumber\\
 \sup_{t\ge 0} \frac{(1+|t|)^{1-2/p}}
 {\log(2+|t|)}\|f(t)\|_{L^p}<\infty,
 \sup_{t\ge 0} \|f(t)\|_{L^2}<\infty\}\nonumber\\
 \nonumber
 \end{eqnarray} endowed with the norm
$$
 \|f\|_Y=\max\{\sup_{t\ge
 0}{(1+|t|)^{1-\frac{2}{p}}}\|f(t)\|_{L^2_{-\sigma}},
 \sup_{t\ge 0} \frac{(1+|t|)^{1-2/p}}
 {\log(2+|t|)}\|f(t)\|_{L^p},
 \sup_{t\ge 0} \|f(t)\|_{L^2}\}
$$

 \noindent To this extent we  consider the operator $N$ defined
on functions in $Y$ as
\begin{eqnarray}
 (Nu)(t)=\Omega(t,s)v+\int_0^t \Omega(t,\tau)\gamma\{P_c[2\psi_E
 |u|^2+\bar\psi_E u^2+|u|^2 u]\nonumber\\
 -Dh|_{a(\tau)}\langle\psi_0,2\psi_E|u|^2+
\bar\psi_E u^2+|u|^2 u\rangle\}d\tau\nonumber
\end{eqnarray}

\par\noindent We will need some properties of the operator $N$ which
are summarized in

\smallskip
\begin{lemma} We have
\par\noindent{\it (i)} The range of $N$ is $Y$ i.e. $N:Y\to Y$ is
 well defined.
 \par\noindent{\it (ii)} There exists $\tilde C>0$ such that
 $$
 \|Nu_1-Nu_2\|_Y\le \tilde
 C(\|u_1\|_Y+\|u_2\|_Y+\|u_1\|_Y^2+\|u_2\|_Y^2) \|u_1-u_2\|_Y
 $$  In particular $N$ is locally Lipschitz.
 \par Moreover, $\tilde C=\tilde C(C,C_p,C_{p,q_0'})$ where the
 constants $C,C_p,C_{p,q_0'}$ are those from the linear
 estimates for $\Omega(t,s)$ (see Theorems ~\ref{theorem:linear1},
 ~\ref{theorem:linear2} in the next section).
 \label{lemma:N}
 \end{lemma}

 \smallskip\par
 \noindent{\bf Proof of Lemma~\ref{lemma:N}}
   Let us observe that it will suffice to show part {\it (ii)} and then
using the fact that
   $N(0)\equiv 0$ we will have part {\it (i)}. Indeed, part {\it (ii)}
will give us that for $u_1,u_2\in Y$
   we have $Nu_1-Nu_2\in Y$. Taking $u_2\equiv 0$ and since $N(0)\equiv
0$ this will imply
   that $Nu_1\in Y$.

  \par Thus, take $u_1,u_2\in Y$ and consider the difference
$Nu_1-Nu_2$, which  is
 \begin{eqnarray}
  (Nu_1-Nu_2)(t)=\int_0^t\Omega(t,\tau)\gamma\{
 P_c[2\psi_E(|u_1|-|u_2|)(|u_1|+|u_2|)+\bar\psi_E
 (u_1-u_2)(u_1+u_2)\nonumber\\
 +(u_1-u_2)|u_1|^2+(|u_1|-|u_2|)(u_2|u_1|+u_2|u_2|)]
 \nonumber\\
 -Dh|_{a(\tau)}\langle\psi_0,2\psi_E(|u_1|-|u_2|)(|u_1|+|u_2|)+
\bar\psi_E (u_1-u_2)(u_1+u_2)+\nonumber\\
(u_1-u_2)|u_1|^2+(|u_1|-|u_2|)
 (u_2|u_1|+u_2|u_2|)\rangle\}d\tau\nonumber
\end{eqnarray}

\noindent\par{\it The $L^2_{-\sigma}$ estimate} We can work under the less
restrictive hypothesis: $4< p<\infty.$ Let $L^{p'}$ be the dual of $L^p,$ i.e.
$1/p'+1/p=1.$ We have
\begin{eqnarray}
 \|Nu_1-Nu_2\|_{L^2_{-\sigma}}
 \le
\int_0^t\|\Omega(t,\tau)\|_{L^2_\sigma\to L^2_{-\sigma}}\times
\nonumber\\
 \times \|\underbrace{2\psi_E<x>^\sigma(|u_1|-
 |u_2|)(|u_1|+|u_2|)
+\bar\psi_E<x>^\sigma(u_1-u_2)(u_1+u_2)}_A\|_{L^2}d\tau\nonumber\\
 +\int_0^t \|\Omega(t,\tau)\|_{L^{p'}\to L^2_{-\sigma}}\times\nonumber\\
 \times
\|\underbrace{(u_1-u_2)|u_1|^2}_{B_1}+
\underbrace{(|u_1|-|u_2|)(u_2|u_1|+u_2|u_2|)}_{B_2}\|_{L^{p'}}
 d\tau\nonumber\\
 +\int_0^t \|\Omega(t,\tau)\|_{L^2_\sigma\to
 L^2_{-\sigma}}\|Dh|_{a(\tau)}\|_{L^2_\sigma}\times\nonumber\\
 \times\{\underbrace{|\langle\psi_0,2\psi_E(|u_1|-
 |u_2|)(|u_1|+|u_2|)+\bar\psi_E(u_1-u_2)(u_1+u_2)>|}_{F}+\nonumber\\
 +\underbrace{|<\psi_0,
 (u_1-u_2)|u_1|^2+(|u_1|-|u_2|)(u_2|u_1|+u_2|u_2|)\rangle|}_{G}\} d\tau
 \label{eq:L2-sigmaN}
 \end{eqnarray}

 To estimate the term  $A$ we observe that
\begin{equation}\label{eq:Aest}
\|<x>^\sigma\psi_E(|u_1|-|u_2|)(|u_1|+|u_2|)\|_{L^2} \le
\|<x>^\sigma
\psi_E\|_{L^\alpha}\|u_1-u_2\|_{L^p}\||u_1|+|u_2|\|_{L^p}\end{equation}
with $\frac{1}{\alpha}+\frac{2}{p}=\frac{1}{2}$. Then
\begin{eqnarray}
\int_0^t \|\Omega(t,\tau)\|_{L^2_\sigma\to
L^2_{-\sigma}}A(\tau)d\tau\nonumber\nonumber\\ \le\int_0^t \frac{
C}{(1+|t-\tau|)\log^2(2+|t-\tau|)}\cdot 3\|\psi_E<x>^\sigma
|u_1-u_2|(|u_1|+|u_2|)\|_{L^2}d\tau\nonumber\\
\le  3C\tilde  C_1
\int_0^t\frac{\log^2(2+|\tau|)}{(1+|t-\tau|)\log^2(2+|t-\tau|)}
\frac{\||u_1|-|u_2|\|_Y}{(1+|\tau|)^{(1-\frac{2}{p})}}\cdot
\frac{\||u_1|+|u_2|\|_Y}{(1+|\tau|)^{(1-\frac{2}{p})}}\nonumber\\
\le 3C\tilde C_1\tilde C_2(\|u_1\|_Y+\|u_2\|_Y)
\frac{\|u_1-u_2\|_Y}{(1+|t|)\log^2(2+|t|)}\nonumber
\end{eqnarray} where for the first inequality we used
{\it Theorem}~\ref{theorem:linear1},part {\it (i)}. The constants
are given by $\tilde C_1=\sup_{t>0} \|<x>^\sigma\psi_E\|_{L^\alpha}$
 and $\tilde C_2=\sup_{t>0}(1+|t|)\log^2(2+|t|)\int_0^t \frac{\log^2(2+|\tau|)}
{(1+|t-\tau|)\log^2(2+|t-\tau|)}\cdot\frac{d\tau}
{(1+|\tau|)^{2-\frac{4}{p}}}<\infty,$ because $p>4.$

\par To estimate the cubic terms $B_1,\ B_2$ we can not use the term $\psi_E$ as
before, and this is what forces us to work in the $L^p$ space. We have: $$
\|(u_1-u_2)|u_1|^2\|_{L^{p'}}\le
\|u_1-u_2\|_{L^p}\|u_1\|_{L^\alpha}^2 $$ respectively
$$
\|(u_1-u_2)(u_2|u_1|+u_2|u_2|)\|_{L^{p'}}\le \|u_1-u_2\|_{L^p}
\|u_2\|_{L^\alpha}(\|u_1\|_{L^\alpha}+\|u_2\|_{L^\alpha}) $$ with
$\frac{2}{\alpha}+\frac{1}{p}=\frac{1}{p'}$. Since $4\le p$ we have
$2\le\alpha\le p.$ Therefore we can again interpolate:
$$ \|u_i\|_{L^\alpha}\le
\|u_i\|_{L^2}^{1-b}\|u_i\|_{L^p}^b, i=1,2, $$
where $\frac{1}{\alpha}=\frac{1-b}{2}+\frac{b}{p}$. Combining these relations we
obtain for $B_1:$
\begin{equation}
\|(u_1-u_2)|u_1|^2\|_{L^{p'}}\le
\|u_1-u_2\|_{L^p}\|u_1\|_{L^2}^{2(1-b)}\|u_1\|_{L^p}^{2b}
\label{eq:estimateLp}
\end{equation} respectively, for $B_2:$
\begin{equation}
\|(u_1-u_2)(u_2|u_1|+u_2|u_2|)\|_{L^{p'}}\le \|u_1-u_2\|_{L^p}
\|u_2\|_{L^2}^{1-b}\|u_2\|_{L^p}^b(\|u_1\|_{L^2}^{1-b}\|u_1\|_{L^p}^b+
\|u_2\|_{L^2}^{1-b}\|u_2\|_{L^p}^b)
\label{estimate:Gtype}
\end{equation}
with
$$
\frac{2(1-b)}{2}+\frac{2b}{p}+\frac{1}{p}=\frac{1}{p'}.$$
A consequence of this relation and of $p<\infty$ is:
\begin{equation}
(1-\frac{2}{p})(1+2b)=1+2/p>1
\label{rel:bp}
\end{equation}
which will play an essential role in what follows.

Thus, the  estimate for the term containing $B_1+B_2$  is
\begin{eqnarray}
\int_0^t \|\Omega(t,\tau)P_c\|_{L^{p'}\to
L^2_{-\sigma}}\|B_1+B_2\|_{L^{p'}}\le C_p
(\|u_1\|_Y^2+\|u_2\|_Y^2)\|u_1-u_2\|_Y\times\nonumber\\
\times\int_0^t \frac{\log
(2+|\tau|)^{(1+2b)}}{|t-\tau|^{1-\frac{2}{p}}}
\cdot \frac{1}{(1+|\tau|)^{(1-\frac{2}{p})(1+2b)}}d\tau\nonumber\\
\le C_p \tilde C_3(\|u_1\|_Y^2+\|u_2\|_Y^2)
\frac{\|u_1-u_2\|_Y}{(1+|t|)^{1-2/p}}\nonumber
\end{eqnarray} where for the first inequality we used
{\it Theorem} ~\ref{theorem:linear1}, part {\it (ii)}, inequalities
\eqref{eq:estimateLp},
\eqref{estimate:Gtype} and the definition of the norm in $Y.$ For the last inequality we used the fact that
$(1-\frac{2}{p})(1+2b)>1$ (see (\ref{rel:bp})) with
 $\tilde
C_3=\sup_{t>0}(1+|t|)^{1-2/p}\int_0^t
\frac{\log(2+|\tau|)^{(1+2b)}}{|t-\tau|^{1-\frac{2}{p}}}\frac{1}
{(1+|\tau|)^{(1-\frac{2}{p})(1+2b)}}d\tau<\infty$.

\par For estimating the term containing $F$ we have
$$
|F|\le 3\|\psi_0\|_{L^\infty}\|\psi_E\|_{L^\alpha}
\|u_1-u_2\|_{L^p}(\|u_1\|_{L^p}+\|u_2\|_{L^p})$$
$\frac{1}{\alpha}+\frac{2}{p}=1$. Then,
 the term containing $F$ is estimated as
the term containing $A$ with $\tilde C_1$ replaced by $\tilde
C_4=\sup_{t>0}\|Dh|_{a(\tau)}\|_{L^2_\sigma}
\|\psi_0\|_{L^\infty} \|\psi_E\|_{L^\alpha}$.

\par We estimate $G$ as
$$
|G|\le\|\psi_0<x>^\sigma\|_{L^\alpha}\|u_1-u_2\|_{L^2_{-\sigma}}
(\|u_1\|_{L^p}^2+\|u_2\|_{L^p}^2+\|u_1\|_{L^p}\|u_2\|_{L^p})$$
 with $\frac{1}{\alpha}+\frac{1}{2}+\frac{2}{p}=1$. Then the term containing $G$ is
estimated as
\begin{eqnarray}
\int_0^t \|\Omega(t,\tau)\|_{L^2_\sigma\to L^2_{-\sigma}}|G|d\tau
\le 3\int_0^t \frac{C\|\psi_0<x>^\sigma\|_{L^\alpha}} {(1+|t-\tau|)
\log^2(2+|t-\tau|)}\cdot
\frac{\log^2(2+|\tau|)}{(1+|\tau|)^{3-6/p}}d\tau\nonumber\\
\le 3C\tilde C_5(\|u_1\|_Y^2+\|u_2\|_Y^2)
\frac{\|u_1-u_2\|_Y}{(1+|t|)\log^2(2+|t|)}\nonumber
\end{eqnarray} with $\tilde
C_5=\|\psi_0<x>^\sigma\|_{L^\alpha} \sup_{t>0} (1+|t|)\log^2(2+|t|)
\int_0^t \frac{\log^2(2+|\tau|)}{(1+|t-\tau|)\log^2(2+|t-\tau|)
(1+|\tau|)^{3-6/p}}d\tau<\infty$ because $p>3.$

\smallskip
\par\noindent{\it The $L^p$ estimate}: With $p',\ q_0',\ q'$ given by Theorem \ref{theorem:linear2}, we have
\begin{eqnarray}
 \|Nu_1-Nu_2\|_{L^p}
 \le
\int_0^t\|\Omega(t,\tau)\|_{L^2_\sigma\to L^p}\times
\nonumber\\
 \times \|\underbrace{2\psi_E<x>^\sigma(|u_1|-
 |u_2|)(|u_1|+|u_2|)
+\bar\psi_E<x>^\sigma(u_1-u_2)(u_1+u_2)}_{A}\|_{L^2}d\tau\nonumber\\
 +\int_0^t \|\Omega(t,\tau)\|_{L^{p'}\cap L^{q_0'}\cap L^{q'}\to
L^p}\times\nonumber\\
 \times
\|\underbrace{(u_1-u_2)|u_1|^2}_{B_1}+
\underbrace{(|u_1|-|u_2|)(u_2|u_1|+u_2|u_2|)}_{B_2}\|_{L^{p'}\cap L^{q_0'}\cap L^{q'}}
 d\tau\nonumber\\
 +\int_0^t \|\Omega(t,\tau)\|_{L^2_\sigma\to
L^p}\|Dh|_{a(\tau)}\|_{L^2_\sigma}\times\nonumber\\
 \times\{\underbrace{|\langle\psi_0,2\psi_E(|u_1|-
 |u_2|)(|u_1|+|u_2|)+\bar\psi_E(u_1-u_2)(u_1+u_2)\rangle|}_{F}+\nonumber\\
 +\underbrace{|\langle\psi_0,
 (u_1-u_2)|u_1|^2+(|u_1|-|u_2|)(u_2|u_1|+u_2|u_2|)\rangle|}_{G}\} d\tau
 \label{eq:LpN}
 \end{eqnarray}

\par The term $A$ can be treated exactly as before and for the $\Omega$ term we use
{\it Theorem}~\ref{theorem:linear1} part {\it (iii)}. Since $1<p',\ q_0',\ q'\le
4/3,$ we can estimate the $B_1,B_2$ terms in each of the norms $L^{p'},\ L^{q_0},\
L^{q'}$ as we did above for their $L^{p'}$ norm only. For $\Omega$ we use {\it
Theorem} ~\ref{theorem:linear2}, part {\it (iii)}. The terms $F$ and $G$ are also
treated as in the previous case. The convolution integrals in \eqref{eq:LpN} will
all decay like $(1+|t|)^{-(1-2/p)}$ except the second one which will have a
logarithmic correction dominated by $\log(2+|t|).$

\medskip
\par\noindent{\it The $L^2$ estimate}: We have
\begin{eqnarray}
 \|Nu_1-Nu_2\|_{L^2}
 \le
\int_0^t\|\Omega(t,\tau)\|_{L^2_\sigma\to L^2}\times
\nonumber\\
 \times \|\underbrace{2\psi_E<x>^\sigma(|u_1|-
 |u_2|)(|u_1|+|u_2|)
+\bar\psi_E<x>^\sigma(u_1-u_2)(u_1+u_2)}_{A}\|_{L^2}d\tau\nonumber\\
 +\int_0^t \|\Omega(t,\tau)\|_{L^{p'}\cap L^2\to
L^2}\times\nonumber\\
 \times
\|\underbrace{(u_1-u_2)|u_1|^2}_{B_1}+
\underbrace{(|u_1|-|u_2|)(u_2|u_1|+u_2|u_2|)}_{B_2}\|_{L^{p'}\cap L^2}
 d\tau\nonumber\\
 +\int_0^t \|\Omega(t,\tau)\|_{L^2_\sigma\to
L^2}\|Dh|_{a(\tau)}\|_{L^2_\sigma}\times\nonumber\\
 \times\{\underbrace{|\langle\psi_0,2\psi_E(|u_1|-
 |u_2|)(|u_1|+|u_2|)+\bar\psi_E(u_1-u_2)(u_1+u_2)\rangle|}_{F}+\nonumber\\
 +\underbrace{|\langle\psi_0,
 (u_1-u_2)|u_1|^2+(|u_1|-|u_2|)(u_2|u_1|+u_2|u_2|)\rangle|}_{G}\} d\tau
 \label{eq:L2N}
 \end{eqnarray}

\par We estimate the term $A$ as in \eqref{eq:Aest} while the estimates in
$L^{p'}$ for $B_1,B_2$ are as in (\ref{eq:estimateLp}) and (\ref{estimate:Gtype}).
For their estimate in $L^2$ norm we use
$$
\|(u_1-u_2)|u_1|^2\|_{L^2}\le
\|u_1-u_2\|_{L^p}\|u_1\|_{L^\alpha}^2 $$ respectively
$$
\|(u_1-u_2)(u_2|u_1|+u_2|u_2|)\|_{L^2}\le \|u_1-u_2\|_{L^p}
\|u_2\|_{L^\alpha}(\|u_1\|_{L^\alpha}+\|u_2\|_{L^\alpha}) $$ with
$\frac{2}{\alpha}+\frac{1}{p}=\frac{1}{2}$. Since $6\le p$ we have
$4\le\alpha\le p.$ Therefore we can again interpolate:
$$ \|u_i\|_{L^\alpha}\le
\|u_i\|_{L^2}^{1-b}\|u_i\|_{L^p}^b, i=1,2, $$
where $\frac{1}{\alpha}=\frac{1-b}{2}+\frac{b}{p}$. Combining these relations we
obtain for $B_1:$
$$
\|(u_1-u_2)|u_1|^2\|_{L^2}\le
\|u_1-u_2\|_{L^p}\|u_1\|_{L^2}^{2(1-b)}\|u_1\|_{L^p}^{2b}
$$ respectively, for $B_2:$
$$
\|(u_1-u_2)(u_2|u_1|+u_2|u_2|)\|_{L^2}\le \|u_1-u_2\|_{L^p}
\|u_2\|_{L^2}^{1-b}\|u_2\|_{L^p}^b(\|u_1\|_{L^2}^{1-b}\|u_1\|_{L^p}^b+
\|u_2\|_{L^2}^{1-b}\|u_2\|_{L^p}^b)
$$
with
$$
\frac{2(1-b)}{2}+\frac{2b}{p}+\frac{1}{p}=\frac{1}{2}.$$
A consequence of this relation is:
$$
(1-\frac{2}{p})(1+2b)=2
$$
Using now the definition of the norm in $Y$ we will have:
$$\|B_1+B_2\|_L^2\le
\|u_1-u_2\|_Y(\|u_1\|_Y^2+\|u\|_2^2)\frac{\log^{2p/(p-2)}(2+|t|)}{(1+|t|)^2}$$

The previous estimates for $F$ and $G$ suffice here as well.

\par Recalling from {\it Theorem }~\ref{theorem:linear1}, part
{\it(iii)} and {\it Theorem }~\ref{theorem:linear2}, part {\it(ii)}, that
$\|\Omega(t,\tau)\|_{L^2_\sigma\to L^2}$ and $\|\Omega(t,\tau)\|_{L^{p'}\cap
L^2\to L^2}$ are bounded, and combining with the estimates above, as well as
taking into account the definition of the functional space $Y$ we have that
$$
\|Nu_1-Nu_2\|_{L^2}\le C\|u_1-u_2\|_Y [\tilde
C_6(\|u_1\|_Y+\|u_2\|_Y)+\tilde C_7(\|u_1\|_Y^2+ \|u_2\|_Y^2)] $$
with $\tilde C_6=\sup_{t\ge 0}\int_0^t
\frac{\log^2(2+|\tau|)}{(1+|\tau|)^{(2-4/p)}}d\tau<\infty$ and $
\tilde C_7=\sup_{t\ge 0}\int_0^t
\frac{\log^{(1+2b)}(2+|\tau|)}{(1+|\tau|)^{(1+2b)(1-2/p)}}d\tau<\infty$.
\par This finishes the proof of {\it Lemma}~\ref{lemma:N}. $\Box$

\smallskip\par We can continue now with the proof of {\it
Lemma}~\ref{lemma:r}. Let
$$
b=\min\{1-\varepsilon,\frac{1-\epsilon}{4\tilde C}\}$$ where $\varepsilon>0$
arbitrary and $\tilde C$ is given by Lemma \ref{lemma:N}. Consider the ball of
radius $b$ centered at $0$
$$
\mathcal{B}=B(0,b)
$$

\par\noindent Using {\it Lemma}~\ref{lemma:N}, part {\it (ii)} with
$u_2=0$ we have, for $u_1\in\mathcal{B}$
$$
\|Nu_1\|_Y\le \tilde C(\|u_1\|_Y^2+\|u_1\|_Y)\|u_1\|_Y\le \tilde C
2b\cdot b<b$$ which means that the ball $\mathcal{B}$ is invariant
under the action of the operator $N$.

\par Also, using  again {\it Lemma}~\ref{lemma:N}, part {\it (ii)},
we have, for $u_1,u_2\in\mathcal{B}$
\begin{eqnarray}
\|Nu_1-Nu_2\|_Y\le \tilde
C(\|u_1\|_Y^2+\|u_2\|_Y^2+\|u_1\|_Y+\|u_2\|_Y)\|u_1-u_2\|_Y\nonumber\\\le
\tilde C4b\|u_1-u_2\|_Y\le (1-\epsilon)\|u_1-u_2\|_Y\nonumber
\end{eqnarray} which shows that $N:\mathcal{B}\to\mathcal{B}$ is
a contraction. This finishes the proof of {\it Lemma}~\ref{lemma:r}
and of {\it Theorem}~\ref{theorem:main}. $\Box$

\section{Linear Estimates}\label{se:lestimates}

\par Consider the linear Schr\" odinger equation with a potential in
two space dimensions: \[
  \begin{cases}
  i\frac{\partial u}{\partial t}=(-\Delta+V(x))u\\
  u(0)=u_0.
  \end{cases}
 \]
It is known that if $V$ satisfies hypothesis (H1)(i) and (ii) then
the radiative part of the solution, i.e. its projection onto the
continuous spectrum of $H=-\Delta +V,$ satisfies the estimates:
 \begin{equation}\label{Murata}
 \|e^{-iHt}P_c u_0\|_{L^2_{-\sigma}}\le C_M \frac{1}{(1+|t|)\log^2(2+|t|)}\|u_0\|_{L^2_\sigma}
 \end{equation}
for some constant $C_M>0$ independent of $u_0$ and $t\in\mathbb{R},$
see \cite[Theorem 7.6]{mm:ae}, and
 \begin{equation}\label{est:Lp}
 \|e^{-iHt}P_c u_0\|_{L^p}\le  \frac{C_p}{|t|^{1-2/p}}\|u_0\|_{L^{p'}}
 \end{equation}
for some constant $C_p>0$ depending only on $p\ge 2$ and $p'$ given
by $p'^{-1}+p^{-1}=1.$ The case $p=\infty$ in (\ref{est:Lp}) is
proven in \cite{ws:de2}. The conservation of the $L^2$ norm, see
\cite[Corollary 4.3.3]{caz:bk}, gives the $p=2$ case:
$$\|e^{-iHt}P_c u_0\|_{L^2}=
\|u_0\|_{L^{2}}.$$ The general result (\ref{est:Lp}) follows from
Riesz-Thorin interpolation.

\par We would like to extend this estimates to the linearized dynamics
around the center manifold. In other words we consider the linear equation, with initial data at time
  $s$,\begin{eqnarray}
  \begin{cases}
  i\frac{\partial z}{\partial t}&=(-\Delta+V(x))z+\gamma
P_c[2|\psi_E(t)|^2
  z+\psi_E^2(t)\bar z
  +Dh|_{a(t)}(\langle\psi_0,2|\psi_E|^2 z+\psi_E^2\bar z\rangle)]\nonumber\\
  z(s)&=v.
  \end{cases}
  \end{eqnarray}
Note that this is a nonautonomous problem as the bound state
$\psi_E$ around which we linearize may change with time.

  \par By Duhamel's principle we have:
 \begin{eqnarray}
 z(t)=e^{-iH(t-s)}P_c v(s)-i\int_s^t e^{-iH(t-\tau)}\gamma
P_c[2|\psi_E|^2
 z+\psi^2_E \bar z+\nonumber\\
 Dh|_{a(\tau)}(\langle\psi_0,2|\psi_E|^2 z+\psi_E^2\bar z\rangle)]d\tau
 \label{rel:Duhamellin}
 \end{eqnarray}

 \par As in (\ref{def:Omega}) we denote
 \begin{equation}\label{def:Omega1}
 \Omega(t,s)v\stackrel{def}{=}z(t).
 \end{equation}
In the next two theorems we will extend estimates of type
(\ref{Murata})-(\ref{est:Lp}) to the operator $\Omega(t,s)$ relying
on the fact that $\psi_E(t)$ is small. It would be useful to find
sufficient conditions under which our results generalize to large
bound states. Such conditions have been obtained in one or three
space dimensions, see \cite{bp:asi, ks:sm1, ws:sm3, sc:as},
unfortunately their techniques cannot be applied in the two space
dimension case.

We start with estimates in weighted $L^2$ spaces:

\begin{theorem} \label{theorem:linear1} There exists $\varepsilon_1>0$ such that if
$\|<x>^\sigma\psi_E\|_{H^2}<\varepsilon_1$ then there exist
constants $C,\ C_p>0$ with the property that for any $t,\
s\in\mathbb{R}$ the following hold: $$ \textrm{(i)}\
\|\Omega(t,s)\|_{L^2_\sigma\to L^2_{-\sigma}}\le
\frac{C}{(1+|t-s|)\log^2(2+|t-s|)}$$

$$\textrm{(ii)}\ \|\Omega(t,s)\|_{L^{p'}\to
L^2_{-\sigma}}\le \frac{C_p}{|t-s|^{1-\frac{2}{p}}},\ {\rm for\
any}\ \infty>p\ge 2\ {\rm where}\ p'^{-1}+p^{-1}=1$$

$$\textrm{(iii)}\ \|\Omega(t,s)\|_{L^2_\sigma\to L^p}\le
\frac{C_p}{|t-s|^{1-\frac{2}{p}}},\ {\rm for\ any}\ p\ge 2$$
\end{theorem}

\smallskip Before proving the theorem let us remark that {\it (i)} is a generalization of (\ref{Murata})
while {\it (ii)} and {\it (iii)} are a mixture between
(\ref{Murata}) and (\ref{est:Lp}). We have used all these estimates
in the previous section. They are consequences of contraction
principles applied to (\ref{rel:Duhamellin}) and involve estimates
for convolution operators based on (\ref{Murata}) and
(\ref{est:Lp}). It will prove much more difficult to remove the
weights from the estimates {\it (ii)} and {\it (iii)}, see Theorem
\ref{theorem:linear2}.

\smallskip \noindent\par{\bf Proof of Theorem ~\ref{theorem:linear1}}
\smallskip

 Fix $s\in\mathbb{R}.$

\par {\it (i)} By definition (see
(\ref{def:Omega1})), we have $\Omega(t,s)v=z(t)$ where $z(t)$
satisfies equation (\ref{rel:Duhamellin}). We are going to prove the
estimate by showing that the nonlinear equation
(\ref{rel:Duhamellin}) can be solved  via contraction principle
argument in an appropriate  functional space. To this extent let us
consider the functional space $$
 X_1:=\{z\in C(\mathbb{R},L^2_{-\sigma}(\mathbb{R}^2))|\sup_{t\in
 \mathbb{R}}(1+|t-s|)\log^2(2+|t-s|)\|z(t)\|_{L^2_{-\sigma}}<\infty\}
 $$ endowed with the norm $$
 \|z\|_{X_1}:=\sup_{t\in\mathbb{R}}\{(1+|t-s|)\log^2(2+|t-s|)
 \|z(t)\|_{L^2_{-\sigma}}\}<\infty
 $$
Note that the inhomogeneous term in (\ref{rel:Duhamellin}):
$$z_0(t)\stackrel{{\rm def}}{=}e^{-iH(t-s)}P_c v$$ satisfies $z_0\in X_1$ and
\begin{equation}\label{eq:z0}
\|z_0\|_{X_1}\le C_M\|v\|_{L^2_\sigma}
\end{equation}
 because of (\ref{Murata}).

 \par We collect the $z$ dependent part of the right hand side of (\ref{rel:Duhamellin}) in
 a linear operator $L(s):X_1\rightarrow X_1 ,$
 $$
 [L(s)z](t)=-i\int_s^t e^{-iH(t-\tau)}\gamma P_c[2|\psi_E|^2
 z+\psi_E^2\bar z+Dh|_{a(\tau)}
(\langle\psi_0,2|\psi_E|^2 z+\psi_E^2\bar z\rangle)]d\tau
$$

  In what
follows we will show that $L$ is a well defined bounded operator from $X_1$ to
$X_1$
 whose operator norm can be made less or equal to $1/2$ by choosing $\varepsilon_1$ in the
 hypothesis sufficiently small. Consequently $Id-L$ is invertible and the solution of the equation
 (\ref{rel:Duhamellin}) can be written as $z=(Id-L)^{-1}z_0.$ In particular
 $$\|z\|_{X_1}\le (1-\|L\|)^{-1}\|z_0\|_{X_1}\le 2\|z_0\|_{X_1}$$
which, in combination with the definition of $\Omega,$ the definition of the norm
in $X_1$ and estimate (\ref{eq:z0}), finishes the proof of {\it (i)}.

It remains to prove that $L$ is a well defined bounded operator from $X_1$ to
$X_1$ whose operator norm can be made less than $1/2$ by choosing $\varepsilon_1$
in the hypothesis sufficiently small. We have the following estimates:
\begin{eqnarray}
 \|L(s) z(t)\|_{L^2_{-\sigma}}\le\int_s^t
 \|e^{-iH(t-\tau)}P_c\|_{L^2_\sigma\to L^2_{-\sigma}}\cdot[3 \|
 |\psi_E|^2(\tau)z(\tau)\|_{L^2_\sigma}\nonumber\\
+\|Dh|_{a(\tau)}\|_{\mathbb{C}\rightarrow
L^2_\sigma}|\langle\psi_0<x>^\sigma,2|\psi_E|^2<x>^{-\sigma}z(\tau)+\psi_E^2
<x>^{-\sigma}\bar z(\tau)\rangle| d\tau\nonumber\\
 \le\int_s^t
 \|e^{-iH(t-\tau)}P_c\|_{L^2_\sigma\to L^2_{-\sigma}}\cdot[3 \|
 |\psi_E|^2(\tau)z(\tau)\|_{L^2_\sigma}\nonumber\\
+\|Dh|_{a(\tau)}\|_{\mathbb{C}\rightarrow
 L^2_\sigma}\|\psi_0\|_{L^2_\sigma}3\|\psi_E^2\|_{L^\infty}
 \|z(\tau)\|_{L^2_{-\sigma}}]\nonumber
 \end{eqnarray}
On the other hand
 \begin{eqnarray}
 \| |\psi_E|^2 z\|_{L^2_\sigma}\le \|z\|_{L^2_{-\sigma}}
 \|<x>^{2\sigma}\|\psi_E|^2\|_{L^\infty},\textrm{ and}\,
\|<x>^\sigma\psi_E\|_{L^\infty}^2\le \varepsilon_1^2
\label{smallnormT}
 \end{eqnarray} where the last inequality holds because of the
 Sobolev imbedding $H^2(\mathbb{R}^2)\subset L^\infty(\mathbb{R}^2)$
  and of the inequality
  $$
 \|<x>^\sigma \psi_E\|_{H^2}\le
 \varepsilon_1.
 $$
Also $$
 \|Dh|_{a(\tau)}\|\le\bar C\varepsilon_1,
\textrm{as}\,|a(\tau)|<\delta.
 $$
Using the last three relations, as well as the estimate (\ref{Murata})
  and the fact that $z\in X_1$ we obtain that
   \begin{eqnarray}
 \|L(s)\|_{X_1\to X_1}\le \varepsilon_1 \sup_{t>0}[
 (1+|t-s|)\log^2(2+|t-s|)\times\nonumber\\\
 \times\underbrace{\int_s^t \frac{1}
 {1+|t-\tau|\log^2(2+|t-\tau|)}\cdot\frac{1}{(1+|\tau-s|)\log^2(2+
 |\tau-s|)}d\tau}_{\mathcal{I}}\le 2C_1\varepsilon_1 \label{normT}
 \end{eqnarray}

 \par Indeed, in order to prove the above we will split $\mathcal{I}$
  into $A+B$ where
  $$
 A=\int_s^{\frac{t+s}{2}}\frac{1}{(1+|t-\tau|)\log^2(2+|t-\tau|)}
 \frac{1}{(1+|\tau-s|)\log^2(2+|\tau-s|)}d\tau
 $$ for which we have the bound
 \begin{eqnarray}
 |A|\le \frac{1}{(1+|\frac{t-s}{2}|)\log^2(2+|\frac{t-s}{2}|)}
 |\int_s^{\frac{t+s}{2}}\frac{1}{(1+|\tau-s|)\log^2(2+|\tau-s|)}|
 \nonumber\\
 \le C_2\frac {1}{(1+|\frac{t-s}{2}|)\log^3(2+|\frac{t-s}{2}|)}
\nonumber
 \end{eqnarray}
Observing that $A=B$ and using the last estimate in
 (\ref{normT}) we obtain that
$$
 \|L\|_{X_1\to X_1}\le C_1\varepsilon_1\le 1/2
$$
for $\varepsilon_1$ small enough.

\smallskip \par {\it (ii)} By the definition of $\Omega$ it is sufficient to prove
that the solution of (\ref{rel:Duhamellin}) satisfies
\begin{equation}\label{WL2-sigmaLp'}
\|z(t)\|_{L^2_{-\sigma}}\le
\frac{C_p}{|t-s|^{1-\frac{2}{p}}}\|v\|_{L^{p'}},\ {\rm for\ all}\ \infty >p\ge 2\ {\rm where}\ p'^{-1}+p^{-1}=1.
\end{equation}
We will use a similar
 functional analytic argument as in the proof of {\it (i)}.
 Fix $p, 2\le p<\infty$ and assume $v\in L^{p'},\ p'^{-1}+p^{-1}=1.$ We will work in the following functional
space:
$$
 X_2:=\{z\in C(\mathbb{R},L^2_{-\sigma}(\mathbb{R}^2)|\sup_{t\in
 \mathbb{R}}\|z(t)\|_{L^2_{-\sigma}}|t-s|^{1-\frac{2}{p}}<\infty\}
$$ endowed with the norm
$$
\|z\|_{X_2}:=\sup_{t\in\mathbb{R}}\|z(t)\|_{L^2_{-\sigma}}|t-s|^{1-\frac{2}{p}}<\infty.
$$

Using the fact that  $L^p\hookrightarrow L^2_{-\sigma}$ continuously
and  the estimate (\ref{est:Lp}) we have $e^{-iH(t-s)}P_cv\in X_2.$
In addition, for $L$ defined in the proof of {\it (i)}, we have
\begin{eqnarray}
 \sup_{t>0}|t-s|^{1-\frac{2}{p}}\|L(s)
 z(t)\|_{L^2_{-\sigma}}\le\nonumber\\
 \sup_{t>0}|t-s|^{1-\frac{2}{p}}\int_s^t
 \|e^{-iH(t-\tau)}P_c\|_{L^2_\sigma\to
  L^2_{-\sigma}}\cdot [\| |\psi_E|^2(\tau)z(\tau)+\nonumber\\
+\psi_E^2(\tau)\bar
  z(\tau)\|_{L^2_\sigma}+\|P_c Dh|_{a(\tau)}\|_{\mathbb{C}\to
L^2_\sigma}|\langle\psi_0,
2|\psi_E|^2z(s)+\psi_E^2\bar z(s)\rangle|d\tau\nonumber\\
  \le\sup_{t>0} |t-s|^{1-\frac{2}{p}}\int_s^t
\frac{(C+\|\psi_0\|_{L^2})\|\psi_E^2<x>^{2\sigma}\|_{L^\infty}}{(1+|t-\tau|)
|\tau-s|^{1-\frac{2}{p}}\log^2(2+|t-\tau|)}<C_3\varepsilon_1^2.\label{normTtris}
  \end{eqnarray}
Using now the
 bounds (\ref{smallnormT}) in (\ref{normTtris}), for $\varepsilon_1$
small
 enough,
  we obtain that the norm of $L(s)$ is less or equal to $1/2$, i.e. the
operator $Id-L(t,s)$ is invertible, which, as in the proof of {\it (i)}, finishes
the proof of estimate {\it (ii)}.

\smallskip\par{\it (iii)} We already know from part {\it (i)} that equation
(\ref{rel:Duhamellin}) has a unique solution in $L^2_{-\sigma}$ provided $v\in
L^2_\sigma.$ We are going to show that the right hand side of
(\ref{rel:Duhamellin}) is in $L^p.$ Indeed
\begin{equation}\label{frhs}
\|e^{-iH(t-s)}P_cv(s)\|_{L^{p}}\le
\frac{C_p}{|t-s|^{1-\frac{2}{p}}}\|v(s)\|_{L^{p'}}\le
\frac{C_p}{|t-s|^{1-\frac{2}{p}}}\|v(s)\|_{L^2_\sigma}
\end{equation} where the $C_p$'s in the two inequalities are
different, for the first inequality we used (\ref{est:Lp}) while for
the second we used the continuous embedding
$L^2_\sigma\hookrightarrow L^{p'},1\le p'\le 2$. For the remaining
terms we combine (\ref{est:Lp}) with $\|z\|_{X_1}<\infty$ obtained
in part {\it (i)}:
\begin{eqnarray}
\|i\int_s^t e^{-iH(t-\tau)}P_c (2|\psi_E^2 z(\tau)+\psi_E^2 \bar
z(\tau))d\tau\|_{L^p}\nonumber\\
\le \int_s^t\frac{2C_4}{|t-\tau|^{1-\frac{2}{p}}}\|<x>^\sigma
\psi_E^2\|_{\alpha}
\|<x>^{-\sigma}z(\tau)\|_{L^2}\nonumber\\
\le
\int_s^t\frac{2C_4}{|t-\tau|^{1-\frac{2}{p}}}\frac{C\varepsilon_1^2}{(1+|\tau-s|)\log^2(2+
|\tau-s|)}d\tau\le \frac{C}{|t-s|^{1-\frac{2}{p}}}\label{srhs}
\end{eqnarray} with $\frac{1}{\alpha}+\frac{1}{2}=\frac{1}{p'}$.

\par Similarly, we have
\begin{eqnarray}
\|i\int_s^t e^{-iH(t-\tau)}\gamma
P_cDh|_{a(\tau)}\langle\psi_0,2\|\psi_E|^2z(s)+\psi_E^2
\bar z(s)\rangle d\tau\|_{L^p}\nonumber\\
\le
\int_s^t\frac{C_5}{|t-\tau|^{1-\frac{2}{p}}}|\langle\psi_0,2\|\psi_E|^2z(s)+\psi_E^2\bar
z(s)\rangle|d\tau\nonumber\\
\le \int_s^t
\frac{C_5}{|t-\tau|^{1-\frac{2}{p}}}|\|\psi_0\|_{L^2}\|<x>^\sigma\psi_E^2\|_
{L^\infty}\|<x>^{-\sigma}z\|_{L^2}\nonumber\\
\le C_6\int_s^t
\frac{1}{|t-\tau|^{1-\frac{2}{p}}}\frac{1}{(1+|\tau-s|)(\log^2(2+|\tau-s|)}d\tau\le
\frac{C}{|t-s|^{1-\frac{2}{p}}}.\label{trhs}
\end{eqnarray}

Plugging (\ref{frhs})-(\ref{trhs}) into (\ref{rel:Duhamellin}) we get:
$$\|z(t)\|_{L^p}\le
\frac{C_p}{|t-s|^{1-\frac{2}{p}}}\|v\|_{L^2_{\sigma}}$$
which by the definition $\Omega(t,s)=z(t)$ finishes the proof of part {\it (iii)}.
$\Box$

The next step is to obtain estimates for $\Omega(t,s)$ in unweighted $L^p$ spaces.
They are needed for controlling the cubic term in the operator $N$ of the previous
section.

\begin{theorem} Assume that $\|<x>^\sigma\psi_E\|_{H^2}<\varepsilon_1$
(where $\varepsilon_1$ is the one used in Theorem~\ref{theorem:linear1}). Then for
all $t,\ s\in\mathbb{R}$ the following estimates hold:
$$
\textrm{(i)}\,\|\Omega(t,s)\|_{L^1\cap L^{q'}\cap L^{p'}\to
L^p}\le \frac{C_{p,q'}\log(2+|t-s|)} {(1+|t-s|)^{1-\frac{2}{p}}},$$ for all $p,\
q',\ 2\le p<\infty,\ 1<q'\le 2,\ p'^{-1}+p^{-1}=1;$

$$
\textrm{(ii)}\,\|\Omega(t,s)\|_{L^2\cap L^{q_0'}\to L^2}\le
C_{q_0'},\ {\rm for\ all}\ q_0',\ 1<q_0'<\frac{4}{3};
$$

\smallskip\par {\it(iii)} for fixed $p_0>0$ and $1<q'_0<4/3$ and for any $2\le p\le
p_0$
$$
\|\Omega(t,s)\|_{L^{q'}\cap L^{p'}\cap L^{q_0'}\to L^p}\le
\frac{C_{p,q_0'}\log(2+|t-s|)^{\frac{1-2/p}{1-2/p_0}}}{|t-s|^{1-\frac{2}{p}}}
$$ where
$$
\frac{1}{q'}=\theta+\frac{1-\theta}{q_0'}
$$
with
$$
\theta=\frac{1-2/p}{1-2/p_0},\ {\textrm i.e.}\ \frac{1}{p}=\frac{\theta}{p_0}+\frac{1-\theta}{2}.
$$
\label{theorem:linear2}
\end{theorem}

Note that {\it (iii)} is similar to the standard estimate for Schr\" odinger
operators (\ref{est:Lp}) except for the logarithmic correction and a smaller
domain of definition. We will obtain it by interpolation from {\it (i)} and {\it
(ii)}. The proof of {\it (i)} will rely on a fixed point technique for equation
(\ref{def:W}) while the proof of {\it (ii)} will rely on Strichartz inequalities.

It turns out that we need to regularize (\ref{rel:Duhamellin}) in
order to obtain {\it (i)} and {\it (ii)}. The inhomogeneous term has
a nonintegrable singularity at $t=s$ when estimated in $L^\infty:$
$$\|e^{-iH(t-s)}P_cv\|_{L^\infty}\le |t-s|^{-1}\|v\|_{L^1}.$$
Using estimates with integrable singularities at $t=s,$ for example
in $L^p,\ p<\infty$ see (\ref{est:Lp}), would lead to a slower time
decay in {\it (i)} and eventually will make it impossible to close
the estimates for the operator $N$ in the previous section. We avoid
this by defining:
\begin{equation}
 W(t)\stackrel{def}{=}[z(t)-e^{-iH(t-s)}P_c]v(s)\label{rel:z}
 \end{equation} which, by plugging in (\ref{rel:Duhamellin}), will satisfy the
following "regularized" equation:
\begin{eqnarray}
 W(t)=\underbrace{-i\int_s^t e^{-iH(t-\tau)}\gamma
P_c[2|\psi_E(\tau)|^2
 e^{-iH(\tau-s)}P_cv(s)+\psi^2_E(\tau)e^{iH(\tau-s)}P_c\bar
 v(s)]d\tau}_{f(t)}\nonumber\\
-\underbrace{i\int_s^t e^{-iH(t-\tau)}\gamma P_cDh|_{a(\tau)}
\langle\psi_0,2|\psi_E|^2e^{-iH(\tau-s)}P_c v(s)+\psi_E^2
e^{iH(\tau-s)}P_c
\bar v(s)\rangle d\tau}_{\tilde f(t)} \nonumber\\
 -\underbrace{i\int_s^t
e^{-iH(t-\tau)}P_c(2|\psi_E|^2W(\tau)+\psi_E^2\bar
 W(\tau))d\tau}_{g(t)}\nonumber\\
-\underbrace{i\int_s^t e^{-iH(t-\tau)} \gamma
P_cDh|_{a(\tau)}\langle\psi_0,2|\psi_E|^2 W(s)+\psi_E^2 \bar W(s)\rangle
d\tau}_{\tilde g(t)} \label{def:W}
 \end{eqnarray}

Some other new notations are necessary for the sake of easy reference. We will
denote by $T(t,s)$ the operator which associates to the initial data at time $s$,
$v$, the function $W(t)$, so that
\begin{equation}\label{def:T}
T(t,s)v\stackrel{def}{=}W(t)
\end{equation}
which will be related to the operator $\Omega(t,s)=z(t)$ (see (\ref{def:Omega}))
by
\begin{equation}\label{rel:TO}
\Omega(t,s)=T(t,s)+e^{-iH(t-s)}P_c.
\end{equation}
For $T$ we can not only extend the estimate in
Theorem~\ref{theorem:linear1} {\it (ii)} to the case $p=\infty$ but
also obtain a nonsingular version of it:

 \begin{lemma}\label{le:T} Assume that $\|<x>^\sigma\psi_E\|_{H^2}<\varepsilon_1$
(where $\varepsilon_1$ is the one used in Theorem~\ref{theorem:linear1}). Then for
each $1< q'\le 2$ there exists the constant $C_{q'}>0,\ C_{q'}\to\infty$ as $q'\to
1,$ such that for all $t,\ s\in\mathbb{R}$
 we have:
$$
 \|T(t,s)\|_{L^1\cap L^{q'}\rightarrow L^2_{-\sigma}}\le \frac{C_{q'}}{1+|t-s|}.
$$
 \end{lemma}

 \smallskip\par {\bf Proof of the Lemma:} Fix $q',\ 1<q'\le 2.$ Consider equation (\ref{def:W})
 with $s\in\mathbb{R}$ arbitrary and $v\in L^1\cap L^{q'}.$ We are going to show
 that (\ref{def:W}) has a unique solution in $C(\mathbb{R},L^2_{-\sigma})$
 satisfying:
 $$\|W(t)\|_{L^2_{-\sigma}}\le\frac{C_{q'}}{1+|t-s|}\max\{\|v\|_{L^1},\
 \|v\|_{L^{q'}}\}$$ which will be equivalent to the conclusion of the Lemma via the
 definition of $T$ (\ref{def:T}).

 Let us observe
 that it suffices to
 prove this estimate only for the forcing term $f(t)+\tilde f(t)$
because then
 we will be able to do the contraction principle in the functional
 space (in time and space) in which $f(t)+\tilde f(t)$ will be, and
thus obtain
 the same decay for $W$ as for $f(t)+\tilde f(t)$.

 \par Indeed, this time  we will consider the functional space
$$
 X_3:=\{u\in C(\mathbb{R},L^2_{-\sigma}(\mathbb{R}^2)|\sup_{t\in
 \mathbb{R}}\|u(t)\|_{L^2_{-\sigma}}(1+|t-s|)<\infty\}
$$ endowed with the norm
$$
\|u\|_{X_3}:=\sup_{t\in\mathbb{R}}\{\|u(t)\|_{L^2_{-\sigma}}(1+|t-s|)\}<\infty
$$

 \par  We have
 \begin{eqnarray}
 \sup_{t>0}(1+|t-s|)\|L(s) u(t)\|_{L^2_{-\sigma}}\le\nonumber\\
 \sup_{t>0}(1+|t-s|)\int_s^t \|e^{-iH(t-\tau)}P_c\|_{L^2_\sigma\to
 L^2_{-\sigma}}
 \times [2\| <x>^\sigma|\psi_E|^2(\tau)<x>^\sigma<x>^{-\sigma}
 u(\tau)\|_{L^2_\sigma}\nonumber\\
+\|Dh\|_{L^2_\sigma}\|\psi_0\|_{L^2}
 \|<x>^\sigma
\psi_E^2\|_{L^\infty}\|u(\tau-s)\|_{L^2_{-\sigma}}]d\tau\nonumber\\
 \le\sup_{t>0} (1+|t-s|)\int_s^t
 \frac{C_7\|\psi_E^2<x>^{2\sigma}\|_{L^\infty}}{(1+|t-\tau|)
 (1+|\tau-s|)\log^2(2+|\tau-s|)}<C_8\varepsilon_1\label{normTbis}
 \end{eqnarray}

 \par Using the
 bounds (\ref{smallnormT}) in (\ref{normTbis}) we obtain that for
 $\varepsilon_1$
 small enough the norm of $L(t,s)$ in $X$ is less then one, i.e. the
operator
 $Id-L(t,s)$ is invertible.

 \par We need now to estimate $f(t)+\tilde f(t)$:
 \begin{equation}
 ||f(t)+\tilde f(t)||_{L^2_{-\sigma}}\le \int_s^t
\frac{2C_M|||\psi_E|^2<x>^\sigma
 e^{-iH(\tau-s)}P_c v(s)||_{L^2}}
 {(1+|t-\tau|)\log^2(2+|t-\tau|)}d\tau \label{rel:fL2sigma}
 \end{equation} (where we used the estimate (\ref{Murata}))

 \par We will split now (\ref{rel:fL2sigma}) into two parts to be
 estimated differently:
 \begin{equation}
 ||f+\tilde f||_{L^2_{-\sigma}}\le
 \underbrace{\int_s^{s+1}\dots}_{\mathcal{I}}+\underbrace{\int_{s+1}^t
 \dots}_{\mathcal{II}}\label{rel:splitf}
 \end{equation}

 \par Then, we have:
 \begin{eqnarray}
 |\mathcal{I}|\le \int_s^{s+1} \frac{2C_M|| |\psi_E|^2<x>^\sigma
 e^{-iH(\tau-s)}P_cv(s)||_{L^2}}{(1+|t-\tau|)\log^2(2+|t-\tau|)}\le
 \nonumber\\
 \le \frac{2C_M}{(1+|t-s-1|)\log^2(2+|t-s-1|)}\int_s^{s+1}
 ||e^{-iH(\tau-s)}P_c v(s)||_{L^q} \cdot \underbrace{||
 |\psi_E|^2<x>^\sigma||_{L^\alpha}}_{\le
 \textrm{fixed constant}}\le\nonumber\\
 \le \frac{C_9}{(1+|t-s-1|)\log^2(2+|t-s-1|)}\int_s^{s+1}
 ||v(s)||_{L^{q'}} (\frac{1}{\tau-s})^{1-\frac{2}{q}} d\tau
 \nonumber\\
 \le \frac{C_{10}||v(s)||_{L^{q'}}}{(1+|t-s-1|)\log^2(2+|t-s-1|)}\le
 C_{11}\frac{1}{1+|t-s|}\|v(s)\|_{L^{q'}}\nonumber
 \end{eqnarray} with $\frac{1}{\alpha}+\frac{1}{q}=\frac{1}{2}$
 and $\frac{1}{q}+\frac{1}{q'}=1$.
 \par For the second integral we have:
 \begin{eqnarray}
 |\mathcal{II}|\le \int_{s+1}^t
 \frac{2C_M|||\psi_E|^2<x>^\sigma||_{L^2}||e^{-iH(\tau-s)}P_c
 v(s)||_{L^\infty}}{(1+|t-\tau|)\log^2(2+|t-\tau|)}d\tau\le\nonumber\\
 \le
\int_{s+1}^t\frac{2C_M|||\psi_E|^2<x>^\sigma||_{L^2}}{(1+|t-\tau|)\log^2(2+|t-\tau|)}\cdot
 \frac{1}{|\tau-s|}||v(s)||_{L^1}d\tau\nonumber\\
 \le \frac{C_{12}}{1+|t-s|}||v(s)||_{L^1}\nonumber
 \end{eqnarray}

 \par Let us observe that the last two estimates are for the case
 when $t>s+1$. If $s<t<s+1$ we have
 \begin{eqnarray}
 ||f+\tilde f||_{L^2_{-\sigma}}\le \int_s^t
\frac{2C_M|||\psi_E|^2<x>^\sigma
 e^{-iH(\tau-s)}P_c v(s)||_{L^2}}
 {(1+|t-\tau|)\log^2(2+|t-\tau|)}d\tau \nonumber\\
 \le 2C_M\int_s^t
 \frac{\|\psi_E^2<x>^\sigma\|_{L^\alpha}\|e^{-iH(\tau-s)}
 P_cv(s)\|_{L^q}}{(1+|t-\tau|)(\log^2(2+|t-\tau|))}d\tau\nonumber\\
 \le C_{13}\int_s^t (\frac{1}{\tau-s})^{1-\frac{2}{q}}d\tau
 \|v(s)\|_{q'}\le C\|v(s)\|_{L^{q'}}\nonumber
 \end{eqnarray} with $\frac{1}{\alpha}+\frac{1}{q}=\frac{1}{2}$.

 \par Combining the last three estimates we get the lemma. $\Box$

 \smallskip
 \par We can now proceed with the proof of Theorem~\ref{theorem:linear2}.
\smallskip

\noindent\par{\bf Proof of Theorem~\ref{theorem:linear2}:}

 {\it (i)} Because of estimate (\ref{est:Lp}) and relation (\ref{rel:TO}) it suffices to prove {\it (i)} for $T(t,s).$

 Consider equation (\ref{def:W}) with arbitrary $s\in\mathbb{R}$ and $v\in L^1\cap
 L^{q'}.$ In the previous Lemma we showed that the solution $W(t)\in
 L^2_{-\sigma}.$ Now we show that it is actually in $L^p$ for all $2\le p<\infty.$
 Fix such a $p.$ Then:
 \begin{eqnarray}
 ||W(t)||_{L^p}\le ||f(t)+\tilde f(t)||_{L^p}+\int_s^t
 ||e^{-iH(t-\tau)}P_c||_{L^{p'}\rightarrow L^p}
(2+\|Dh|_{a(\tau)}\|_{L^p}\|\psi_0\|_{L^p}) \nonumber\\
|||\psi_E|^2 <x>^\sigma <x>^{-\sigma}|W(\tau)|||_{L^{p'}}
d\tau\nonumber\\
 \le ||f(t)+\tilde f(t)||_{L^p}+\int_s^t \frac
{2C_{14}}{(t-\tau)^{1-\frac{2}{p}}}
 |||\psi_E|^2 <x>^\sigma||_{L^\alpha}||<x>^{-\sigma}W(\tau)||_{L^2}
 \label{estimateWL^p}
 \end{eqnarray} (with $\frac{1}{\alpha}+\frac{1}{2}=\frac{1}{p'}$)

 \par The estimate for $f(t)+\tilde f(t)$ is similar , but this time
the term
 $|||\psi_E|^2 e^{-iH(\tau-s)}P_c v||_{L^{p'}}$ is controlled
 for $s+1<\tau<t$  by
 $$||\psi_E^2||_{L^{p'}} ||e^{-iH(\tau-s)}P_c
 v||_\infty\le \frac{C_{15}}{|t-s|}||v||_{L^1}
$$ and for $s\le \tau<s+1$ by

 \begin{displaymath}
 \|\psi_E\|_{L^\alpha} \|e^{-iH(\tau-s)}P_cv\|_{L^{q}}\le
  \frac{C_{16}}{(\tau-s)^{1-\frac{2}{q}}}
 \end{displaymath} where $\alpha^{-1}+q^{-1}=p'^{-1}$ and $q^{-1}+q'^{-1}=1.$

 \par Using now the previous Lemma to estimate the
 term $||<x>^{-\sigma}W(\tau)||_{L^2}$ and replacing in
 (\ref{estimateWL^p}) we get:
 $$
 ||W(t)||_{L^p}\le
\frac{C_{17}\log(1+|t-s|)}{(1+|t-s|)^{1-\frac{2}{p}}}
 \max\{||v||_{L^1},||v||_{L^{q'}}\}
 $$ with $1<q'\le 2$ which is equivalent to
\begin{equation}\label{est:Tp}
\|T(t,s)\|_{L^1\cap L^{q'}\rightarrow L^p}\le \frac{C_{17}\log(1+|t-s|)}{(1+|t-s|)^{1-\frac{2}{p}}}
\end{equation}
for all $2\le p<\infty$ and $1<q'\le 2.$ This finishes the proof for
part {\it (i).}

 \smallskip\par {\it (ii)}
  Recalling the equation for $W$
 (\ref{def:W}), let us observe that we have
 \begin{eqnarray}
 \|\int_s^t e^{-iH(t-\tau)}P_c(|\psi_E|^2W(\tau)+\psi_E^2 \bar
 W(\tau))d\tau\|_{L^2}\le
 C_S(\int_s^t
 \|\psi_E|^2W(\tau)\|_{L^{\rho'}}^{\gamma'})^{\frac{1}{\gamma'}}
 \nonumber\\
 \le C_S(\int_s^t
\|\psi_E^2<x>^\sigma\|_{L^{\frac{2\gamma'}{2-\gamma'}}}^
 {\gamma'}
 \|W(\tau)<x>^{-\sigma}\|_{L^2}^{\gamma'}
 d\tau)^{\frac{1}{\gamma'}}\le C_S\varepsilon_1^2\int_s^t
\frac{1}{(1+|\tau-s|)}^
 {\gamma'(1-\frac{2}{q_0})}d\tau\nonumber\\
 \le C_{18}\|v\|_{q_0'}\label{estL21}
 \end{eqnarray} where for the first inequality we used the
 Strichartz  estimate
 $$
 (\mathcal{T}f)(t)=\int_s^t e^{-iH(t-\tau)}f(\tau)d\tau :
   L^{\gamma'}(0,T;L^{\rho'})\to L^{\infty}(0,T;L^2)
 $$ with $(\gamma,\rho)$ with $\gamma\ge 2$ . For the second inequality we
used  H\"{o}lder's inequality and for the third one we used
(\ref{WL2-sigmaLp'}) combined with (\ref{rel:z}) and (\ref{est:Lp}).
    Finally the last inequality
   holds when $\gamma'(1-\frac{2}{q_0})>1$
   which happens for $q_0>2\gamma\ge 4$.

\par Also, we have the estimates
 \begin{eqnarray}
  \|\int_s^te^{-iH(t-\tau)}\gamma P_c
Dh|_{a(\tau)}\langle\psi_0,2|\psi_E|^2W(\tau)+\psi_E^2
 \bar W(\tau)\rangle d\tau\|_{L^2}\nonumber\\
 \le C_S(\int_s^t \|Dh|_{a(\tau)}\langle\psi_0,2|\psi_E|^2W(\tau)+\psi_E^2
 \bar W(\tau)\rangle\|_{L^{\gamma'}}^{\rho'}d\tau)^{
 \frac{1}{\gamma'}}\nonumber\\
 \le \varepsilon_1 C_{19}(\int_s^t
|\langle\psi_0,2|\psi_E|^2W(\tau)+\psi_E^2\bar W(\tau)\rangle|^{
 \gamma'})^{\frac{1}{\gamma'}}\nonumber\\ \le
 \varepsilon_1 C_{20}(\int_s^t
(\|\psi_0\|_{L^2}\|\psi_E^2<x>^\sigma\|_{L^\infty}
 \|W(\tau)\|_{L^2_{-\sigma}})^{\gamma'}d\tau)^{\frac{1}{\gamma'}}
 \nonumber\\
\le \varepsilon_1 C_{21}(\int_s^t
\frac{1}{(1+|\tau-s|)^{(1-\frac{2}{q_0})\gamma'}}
 d\tau)^{\frac{1}{\gamma'}}\|v(s)\|_{L^{q_0'}}\le \varepsilon_1
C_{22}\|v\|_{L^{q_0'}}\label{estL22}
\end{eqnarray} where for the first inequality we used Strichartz
estimate as before and for the second inequality we use the fact
that $Dh|_{a(\tau)}$ is bounded in $H^2$ and thus in any $L^p$ and
its norm is small. For the fourth inequality we used the fact
 that $\|\psi_0\|_{L^2}$ and $\||\psi_E|^2<x>^\sigma\|_{L^\infty}$ are
bounded and small.
 Finally the last inequality holds, as before, for $q_0>2\gamma\ge 4$.

 For $f(t)+\tilde f(t)$ we'll need to estimate differently the
short time
 behavior and the long time behavior, namely:
$$
 f(t)+\tilde
f(t)=\underbrace{\int_s^{s+1}\dots}_{\mathcal{I}}+\underbrace{\int_{s+1}^t
 \dots}_{\mathcal{II}}
$$

 \smallskip\par We have:
 \begin{eqnarray}
 |\mathcal{I}|_{L^2}\le 2\|\int_s^{s+1}
 e^{-iH(t-\tau)}\gamma P_c[|\psi_E|^2e^{-iH(\tau-s)}P_c
 v(s)+Dh(\langle \psi_0,2|\psi_E|^2 e^{-iH(\tau-s)}v(s)\rangle)]\|_{L^2}
 \nonumber\\\le C_{23}\int_s^{s+1} ||
|\psi_E|^2e^{-iH(\tau-s)}P_cv(s)||_{L^2}+|\langle\psi_0,2|\psi_E|^2e^{-iH(\tau-s)}v(s)\rangle|
 d\tau
 \nonumber\\
 \le C_{24}\int_s^{s+1}|||\psi_E|^2||_{L^\alpha}
 ||e^{-iH(\tau-s)}P_c v(s)||_{L^{q_0}}+\|\psi_0\|_{L^2}
 \||\psi_E|^2\|_{L^\alpha}\|e^{-iH(\tau-s)}v(s)\|_{L^{q_0}}d\tau
 \nonumber\\
 \le C_{25}\int_s^{s+1}\frac {1}{(\tau-s)^{1-\frac{2}{q_0}}}d\tau
 ||v||_{L^{q_0'}}\le C_{26}||v||_{L^{q_0'}}\nonumber
 \end{eqnarray} (where we used the fact that the operator $e^{-iHt}$
 preserves the $L^2$ norm, and
 $\frac{1}{\alpha}+\frac{1}{q_0}=\frac{1}{2}$). One estimates
similarly the terms containing $\psi_E^2$ instead of $|\psi_E|^2.$

 We continue by estimating $\mathcal{II}$:
 \begin{eqnarray}
 |\mathcal{II}|_{L^2}\le \|\int_{s+1}^t |\psi_E|^2
 e^{-iH(\tau-s)}P_c
 vd\tau\|_{L^2}\nonumber\\
 \le(\int_{s+1}^t |||\psi_E|^2 e^{-iH(\tau-s)}P_c
 v||_{q_0'}^{\alpha'}d\tau)^{\frac{1}{\alpha'}}\nonumber\\
+(\int_s^t
\|\psi_0\|_{L^2}\||\psi_E|^2\|_{L^\alpha}\|e^{-iH(\tau-s)}
 v\|_{L^{q_0}}^{\alpha'}d\tau)^{\frac{1}{\alpha'}}\nonumber\\
 \le C_{27}(\int_{s+1}^t \|e^{-iH(\tau-s)}P_c
 v\|_{q_0}^{\alpha'})^{\frac{1}{\alpha'}}\le
 C_{28}(\int_{s+1}^t
\frac{1}{(\tau-s)^{\alpha'(1-\frac{2}{q_0})}}d\tau)^{\frac{1}{\alpha'}}
 \|v\|_{q_0'}\le C_{29}\|v\|_{q_0'}\nonumber
 \end{eqnarray} where for the first inequality we used the fact
 that the $L^2$ norm is preserved by the operator $e^{-iHt}P_c$.
 For the second inequality we used the Strichartz estimate
 \begin{displaymath}(\mathcal{T}f)(t)=\int_s^t
e^{-iH(t-\tau)}f(\tau)d\tau
 :L^{q_0'}(0,T;L^{\alpha'})\to L^\infty (0,T;L^2)
 \end{displaymath} for
the $f(t)$ term. For the $\tilde f(t)$ term we used  similarly the same Strichartz
 estimate, the fact that $\|Dh\|_{L^{q_0'}}$ is bounded
(as it is in any $L^p$ norm), and we estimated the scalar product by
the product
 $\|\psi_0\|_{L^2}\||\psi_E|^2\|_{L^\alpha}\|
e^{-iH(\tau-s)} v(s)\|_{L^{q_0}}$. For the
 third inequality we used H\"{o}lder's inequality and the fact that
 the
$\|\psi_E|^2\|_{L^\beta},\||\psi_E|^2\|_{L^\alpha},\|\psi_0\|_{L^2}
 \le C,\forall t$ ( where
 $\frac{1}{q_0'}=\frac{1}{q_0}+\frac{1}{\beta}$ and
$\frac{1}{2}=\frac{1}{\alpha}+ \frac{1}{q_0}$). Finally the last
 inequality holds because $\alpha'(1-\frac{2}{q_0})>1$, as
 $q_0>\frac{2}{\alpha}$.
 \par Let us observe that we assumed that $t>s+1$. If $s<t<s+1$,
 only the estimate for $\mathcal{I}$ will suffice, where the upper
limit of
 integration $s+1$ should be replaced by $t$.

 Combining the estimates for $\mathcal{I},$ $\mathcal{II},$ (\ref{estL21}) and (\ref{estL22}) we
have that $W(t)$ is uniformly bounded in $L^2$ which, by (\ref{def:T}), implies
\begin{equation}\label{est:T2}
\|T(t,s)\|_{L^{q_0'}\rightarrow L^2}\le C_{q_0'},\ {\rm for\ all}\ t,\
s\in\mathbb{R}.
\end{equation}
 Using now (\ref{rel:TO}) and (\ref{est:Lp}) with $p=p'=2$ we obtain {\it (ii)}. $\Box$

\par {\it (iii)} We start from (\ref{est:Tp}):
$$
\|T(t,s)\|_{L^1\cap L^{q_0'}\to L^p}\le
\frac{C_{p,q_0'}\log(2+|t-s|)} {(1+|t-s|)^{1-\frac{2}{p}}} $$
 and
(\ref{est:T2}):
$$
\|T(t,s)\|_{ L^{q_0'}\to L^2}\le C_{q_0'},\ 1<q_0'<\frac{4}{3}.
$$
\par We can now use  the  Riesz-Thorin interpolation between the
spaces $L^1\cap L^{q_0'}$ and $L^{q_0'}$ as starting spaces and
between $L^p$ and $L^2$ as arrival spaces to get the claimed
estimate. Indeed, it suffices to take as in the statement
$\theta=\frac{1-2/p}{1-2/p_0}$, and use it with the above two
relations to get
$$
\|T(t,s)\|_{L^{q'}\cap L^{p'}\cap L^{q_0'}\to L^p}\le
\frac{C_{p,q_0'}\log(2+|t-s|)^{\frac{1-2/p}{1-2/p_0}}}{|t-s|^{1-\frac{2}{p}}}
$$ with $1<q_0'<\frac{4}{3},p\ge 2$ and
$$
\begin{cases}
\frac{1}{q'}=\theta+\frac{1-\theta}{q_0'}\cr
\frac{1}{p}=\frac{\theta}{p_0}+\frac{1-\theta}{2},\theta=\frac{1-2/p}{1-2/p_0}\cr
\end{cases}
$$

\par Using now (\ref{rel:TO}) and (\ref{est:Lp}) we obtain the claimed estimate for $\Omega(t,s)$. $\Box$

\section{Conclusions.}\label{se:conclusions}
%discussion of the evolution on center manifold, large initial data,
%time dependent perturbations, 1-d case.

\par We have established that the solution starting from small and localized
initial data will approach, as $t\to\pm\infty$, the center manifold formed by the
nonlinear bound states (solitary waves). However we have not been able to decide
whether the solution will approach exactly one solitary wave as in the 1-d and 3-d
case, see for example \cite{bs:as,pw:cm}. Here is the main reason:

\par  The long time dynamics on the center manifold is given by the equation
(\ref{eq:at}). Since
$$a(\pm\infty)-a(0)=\lim_{t\rightarrow\pm\infty}\int_0^t\frac{da}{dt}dt,$$ the
existence of an asymptotic limit at $t=\pm\infty$ is equivalent to the
integrability of the right hand side of $\eqref{eq:at}$ on $[0,+\infty)$
respectively $(-\infty,0].$ The terms containing $r^2$ and $r^3$ are absolutely
integrable because they are dominated by  $(1+|t|)^{2(2/p-1)}$, respectively
$(1+|t|)^{3(2/p-1)}$, which are integrable on $\mathbb{R}$ for $p>4$. However, the
linear terms in $r$ do not decay fast enough to be absolutely integrable. It is
possible though that a combination of decay and oscillatory cancellations would
render it integrable. We think that is only a matter of time until a suitable
treatment of this term is found. Note that, in the 1-d and 3-d cases, the linear
terms in $r$ were absolutely integrable in time, see for example
\cite{bs:as,pw:cm}. But these estimates relied on the integrable decay in time of
the Schr\"{o}dinger operator in $L^\infty$ norm in 3-d, respectively on the large
power nonlinearity to compensate for the linear growth in time introduced by
virial type estimates in 1-d. None would work for our cubic NLS in 2-d.

The situation is even more complex and possible more interesting when the center
manifold has more than one branch (more than one connected component). For
simplicity, consider the case when hypothesis (H1) part (iii) is relaxed to allow
for two, simple, negative eigenvalues $E_0<E_1$ with corresponding normalized
eigenvectors $\psi_0,\ \psi_1.$ In this case the center manifold has two branches
$\psi_{E_j}=a_j\psi_j+h_j(a_j),\ j=0,1,$ each bifurcating from one eigenvector as
described in Section \ref{se:prelim}. The decomposition into the evolution on the
center manifold and the one away from it will now be:
$$
u(t,x)=\underbrace{\sum_{j=0}^{1} (a_j(t)\psi_j(x)+h_j(a_j(t)))}_{\psi_{{\rm
CM}}(t)}+r_m(t,x)
$$
The equation for $r_m(t)$ remains essentially the same as (\ref{eq:r1}) in Section
\ref{se:result}, with $\psi_E$ replaced by $\psi_{{\rm CM}}$ and the differential of $h$
replaced by the sum of the differentials of $h_j,\ j=0,1.$ However, one has to add
to the right hand side of \eqref{eq:r1} the projection onto the continuous
spectrum of the interaction term between the branches:
\begin{equation}\label{eq:bit}
2|\psi_{E_0}|^2
\psi_{E_1}+\psi_{E_0}^2\bar\psi_{E_1}+
 2\psi_{E_0}|\psi_{E_1}|^2+\bar\psi_{E_0} \psi_{E_1}^2
\end{equation}
In principle one could use our techniques and obtain a decay in time for $r_m(t),$
hence collapse on the center manifold, provided one makes the ansatz that the term
above, or at least its projection onto the continuous spectrum, decays in time.
Such an ansatz needs to be supported by the analysis of the motion on the center
manifold given now by a system of two ODE's, one for $a_0$ and one for $a_1.$ Each
of the equations will be similar to
\eqref{eq:at} but the projection of \eqref{eq:bit} onto $\psi_0,$ respectively
$\psi_1,$ has to be added to the right hand side. Note that, in the 3-d case,
under the additional assumption $2E_1-E_0>0,$ it has been shown that the evolution
approaches asymptotically a ground state (a periodic solution on the branch
bifurcating from $\psi_0$) except when the initial data is on a finite dimensional
manifold near the excited state branch (the one bifurcating from $\psi_1)$, see
\cite{sw:sgs,ty:ad,ty:rel,ty:sd}. But the authors' analysis relies
heavily on the much better dispersive estimates for Schr\"{o}dinger operators in
3-d compared to 2-d. The 2-d case remains open.

Returning now to the case of one branch center manifold in 2-d, an important
question is whether its stability persists under time dependent perturbations. In
\cite{ckp:res} we showed that this is not the case in 3-d. The slower decay in time of the
Schr\"{o}dinger operator in 2-d compared to 3-d  prevents us, yet again, from
extending the technique in \cite{ckp:res} to the 2-d setting.

\bigskip
\noindent{\bf Acknowledgements:} The authors wish to thank M. I. Weinstein
for helpful comments on the manuscript.  E. Kirr was partially supported by NSF
grants DMS-0405921 and DMS-0603722.

\bibliographystyle{plain}
\bibliography{ref}

\end{document}